\documentclass[11pt]{article}
\usepackage{amssymb,amsmath,verbatim,amsthm}

\newtheorem{theorem}{Theorem}[section]
\newtheorem{lemma}[theorem]{Lemma}
\newtheorem{proposition}[theorem]{Proposition}
\newtheorem{corollary}[theorem]{Corollary}

\newcommand{\ignore}[1]{}

\newcommand{\wdesc}{{\operatorname{wdesc}}}
\newcommand{\desc}{{\operatorname{desc}}}

\begin{document}

\title{Wide-Sense 2-Frameproof Codes\footnote{Supported by NSFC grants 11571034 and 11971053.}}
\author{\small  Junling Zhou, Wenling Zhou\\
\small Department of Mathematics \\
\small Beijing Jiaotong University\\
\small Beijing 100044, P. R. China\\
\small jlzhou@bjtu.edu.cn\\
}

\date{ }
\maketitle

\begin{abstract} Various kinds of fingerprinting codes and their related
combinatorial structures are extensively studied for
protecting copyrighted materials.  This paper concentrates on one specialised fingerprinting code named wide-sense frameproof codes in order to prevent innocent users from being framed.

 Let $Q$ be a finite alphabet of size $q$. Given a $t$-subset $X=\{ x ^1,\ldots, x ^t
\}\subseteq Q^n$, a position $i$ is called
 undetectable for $X$ if the values of the words of $X$ match in
their $i$th position: $x_i^1=\cdots=x_i^t$. The wide-sense descendant set of $X$ is defined
by
$\wdesc(X)=\{y\in Q^n:y_i=x_i^1,i\in {U}(X)\},$  where ${U}(X)$ is the set of undetectable positions for $X$. A code ${\cal C}\subseteq Q^n$ is called  a  wide-sense $t$-frameproof code if
$\wdesc(X) \cap{\cal C} = X$
for all $X \subseteq {\cal C}$ with $|X| \le t$.

The paper improves the upper bounds on the sizes of  wide-sense $2$-frameproof codes by applying techniques  on  non $2$-covering Sperner  families and intersecting families  in extremal set theory.


\end{abstract}


\section{Introduction}

Fingerprinting codes are combinatorial objects that have been studied for more
than 20 years due to their applications in digital data copyright protection and their
combinatorial interest.
 Let $Q$ be a finite alphabet of size $q$. In order to protect  a copyrighted digital product,  a dealer inserts a fingerprint in each copy  and then
distributes  copies to all registered users, where
a {\it fingerprint} is  a string $ x =(x_1,\ldots,x_n)$ over $Q$.  The goal of inserting
the fingerprint is to personalize the copy given out to the user,
and to rule out redistribution. Clearly, an individual user will be deterred from releasing an unauthorized
copy. However, a coalition of  users may collude in order to produce
an unregistered copy.  The
goal of the coalition is to create a fingerprint of
the illegal copy  that is unable to identify users from it.
We assume that the members of a coalition can only alter
those coordinates of the fingerprint in which at least two of their
fingerprints differ, and refer to this as the {\it Marking Assumption}.
 In this paper we concentrate on  $t$-frameproof
codes, which have the property that no coalition of at most $t$ users can frame
a user not in the coalition.

Given a $t$-subset $X=\{ x ^1,\ldots, x ^t
\}\subseteq Q^n$ we now define the the set of descendants of $X$.
   We write $x_i^j$ for the $i$-th component of  $ x ^j$ for $1\le i\le n$ and $1\le j\le t$. A position $i$ is called
{\it undetectable} for $X$ if the values of the words of $X$ match in
their $i$-th position: $x_i^1=\cdots=x_i^t$.
Denote by ${U}(X)$ the set of positions undetectable for $X$. By
the marking assumption, the coalition cannot change the values
of undetectable positions. If the position is detectable, then there
are several options for the coalition to fill it. We will consider
 {\it wide-sense descendant set} defined
by
$$\wdesc(X)=\{y\in Q^n:y_i=x_i^1,i\in {U}(X)\},$$  in contrast to the
{\it  narrow-sense descendant set} $$\desc(X)=\{y\in Q^n:y_i\in\{x_i^1,\ldots,x_i^t\}\}.$$
 In the
literature the wide-sense descendant set could be found under the name of envelope \cite{BargB},
feasible set \cite{BS,StinsonWei}, or  Boneh-Shaw descendant \cite{Blackburnsurvey}.

A code ${\cal C}\subseteq Q^n$ of size $|{\cal C}|=m$ is called an $(n,m,q)$ code. We define ${\cal C}$ to be a {\it wide-sense $t$-frameproof code},
or {\it  $(n,m,q)$ $t$-wFP code}, if
$$\wdesc(X) \cap{\cal C} = X$$
for all $X \subseteq {\cal C}$ with $|X| \le t$. Replacing $\wdesc(X)$ with $\desc(X)$ defines a {\it narrow-sense $t$-frameproof code}.


Boneh and Shaw \cite{BS} were the first to define a $t$-frameproof code, where they adopted the  wide-sense model of descendent sets.
In the narrow-sense model, variants of frameproof codes were extensively studied by many researchers.  Named after the security properties it guarantees, the following types of fingerprinting codes are well-known: frameproof codes,  secure frameproof codes, identifiable parent property codes, traceability codes, and anti-collusion codes, see \cite{IPP}, \cite{BargB}-\cite{{TRACE}}, \cite{CheePhd}-\cite{44}, \cite{FiatTassa,Guo}, \cite{Shangg2}-\cite{Trung}, \cite{Xing-Frame,Yangyt}.

It is clear that $\wdesc(X)$ always strictly contains $\desc(X)$   if $2\le |X|<q$. Just as Blackburn \cite{Blackburnsurvey} said, ``This is one
reason why the problem of constructing analogues of the codes  for
Boneh-Shaw descendants is often more difficult than the original problem." To the best of our  knowledge,  frameproof codes is the only type of fingerprinting codes that was ever studied  in the wide-sense model. 

The paper is organized as follows. In Section 2 we introduce  Sperner families and intersecting families and we display an upper bound on the size of a Sperner family by its maximum size and minimum size of subsets. In Section 3 we improve the known upper bounds on  the sizes of $2$-wFP codes, which were previously established by Panoui in her PhD dissertation \cite{panoui}. Non $2$-covering Sperner families generated by all codewords are considered and better upper bounds are established by developing many results on Sperner families and intersecting families.  In Section 4 we conclude the paper.

\section{Sperner families}

Stinson and Wei \cite{StinsonWei} were the first to establish the relationship between Sperner families and $t$-wFP  codes and then proved that $m\leq {n\choose \lfloor{n\over 2}\rfloor}+1$ for $(n,m,2)$ 2-wFP codes by applying Sperner's Theorem (\cite[Theorem 5.2]{StinsonWei}). Panoui \cite{panoui} developed this idea and presented the equivalence between a $2$-wFP  code and the non $2$-covering  Sperner families generated by all codewords. The upper bounds on the sizes of   $2$-wFP  codes were then improved as follows.

{\lemma\cite[Theorem 6.3.8]{panoui}\label{PwFP} Let ${\cal C}$ be an $(n,m,q)$ $2$-wFP code.

(1) If $n$ is even, then  $m\leq {n\choose {n\over 2}-1}+1$.

(2) If $n$ is odd, then  $m\leq {n\choose {n-1\over 2}}-{n-1\over 2}$.

}

The aim of this section is to  improve the above upper bounds. We first introduce  related definitions in extremal set theory and recall or develop some useful results.

Let ${\cal F}$ be a family of finite sets. If any two distinct sets of ${\cal F}$ are incomparable, that is,  $A\not\subseteq B$ for any $A,B\in{\cal F}$, then ${\cal F}$ is called an {\it antichain} or  a {\it Sperner family}. To the other extreme,  a {\it chain} is a set family ${\cal F}$ in which every pair of sets is comparable.

\begin{theorem}(Sperner's Theorem) \cite{Anderson}
Let ${\cal F}$ be a Sperner family over
an $n$-set. Then  $$|\mathcal{F}| \leq {n \choose {\lfloor\frac{n}{2}\rfloor}}.$$
 \end{theorem}

The size of a Sperner family which contains a singleton is easily obtained from Sperner's Theorem.

{\proposition\cite[Proposition 6.3.4]{panoui}\label{single} Let ${\cal F}$ be a Sperner family over an $n$-set. If there exists
a set $F\in{\cal F}$ such that $|F|= 1$, then $$|{\cal F}|\leq {n-1\choose \lfloor{n-1\over 2}\rfloor}+1.$$
}


Let $C =\{A_1, A_2,\ldots, A_k\}$ be a chain of subsets of  an $n$-set, i.e., $A_1 \subseteq A_2 \subseteq \ldots \subseteq A_k$. This
chain is {\it symmetric} if $|A_1| + |A_k| = n$ and $|A_{i+1}| = |A_i| + 1$ for all $i =
1,2,\ldots, k-1$.

{\theorem \label{symchains}\cite[Theorem 8.3]{Extrem} The family of all subsets of an $n$-set can be partitioned
into ${n \choose {\lfloor\frac{n}{2}}\rfloor}$
mutually disjoint symmetric chains.}

 A family ${\cal F}$ of sets  is called {\it  $k$-intersecting} $(k\ge 1$),  if $|A\cap B|\ge k$  for
all $A,B\in{\cal F}$. An {\it intersecting family} is a 1-intersecting family. Call the families  ${\cal A}$ and ${\cal B}$ {\it cross-$k$-intersecting} if $|A\cap B|\ge k$ holds
for all $A\in {\cal A}$ and $B\in{\cal B}$. ${\cal A}$ and ${\cal B}$ are {\it cross-intersecting} if they are cross-1-intersecting.
 Let ${\cal F}$ be a family of subsets of a ground set $E$. Then ${\cal F}$ is called
{\it non $2$-covering} if for every pair of sets $A,B \in{\cal F}$ we have $A\cup B \ne E$.

{\theorem\label{k-int}
\cite{Milner1968} If  ${\cal F}$ is a $k$-intersecting  Sperner family over an
$n$-set, then $$|{\cal F}|\le{n\choose\lfloor {n+k+1\over 2}\rfloor}.$$
}

Let ${\cal C}= \{c^1,c^2,\ldots,c^m\}$ be an $(n,m,q)$ code. For any $1\le i,j\le m$ and $i\ne j$, define $I(i,j)$ to be the {\it coincidence set} of   $c^i$ and $c^j$, i.e., $$I(i,j)=\{k:c^{i}_k=c^{j}_k,1\le k\le n\}.$$ For any $1\le i\le m$, define $${\cal X}_i=\{I(i,j): i\ne j,1\le j\le m\}$$ to be the {\it coincidence family} generated by the codeword $c^i\in{\cal C}$. Clearly $|{\cal C}|=|{\cal X}_i|+1$ for $1\le i\le m$.

{\theorem\label{FP&Sp} \cite[Lemma 6.3.2, Corollary 6.3.3]{panoui} Let ${\cal C} = \{c^1,c^2,\ldots,c^m\}$ be an $(n,m,q)$ code. Then, ${\cal C}$ is a $2$-wFP code if and only if
 ${\cal X}_i$ is a non $2$-covering Sperner family for any $1\le i\le m$.
 }

%
%
%

%
%

We have a simple but useful result on the coincidence sets.

\begin{lemma}{\label{BH}}
Let ${\cal C}$ be an $(n,m,q)$ code.  For any three codewords $c^{i},c^{j},c^{k}\in \mathcal{C}$, we have
\begin{eqnarray}
 I\left( {i,j} \right) \cap I\left( {i,k} \right) \subseteq  I\left( {j,k} \right) \subseteq  ({I\left( {i,j} \right) \cap I\left( {i,k} \right)})  \cup \overline{I\left( {i,j} \right) \cup I\left( {i,k} \right)}. \notag
\end{eqnarray}
\end{lemma}

\proof Firstly let $p \in  I\left( {i,j} \right) \cap I\left( {i,k} \right) $. Then we have $c_p^i = c_p^j$ and $c_p^i = c_p^k$. Hence
$p \in   I\left( {j,k} \right)$ and $$ I\left( {i,j} \right) \cap I\left( {i,k} \right) \subseteq  I\left( {j,k} \right).$$

Secondly let  $p\in I\left( {j,k} \right)$ and $p\notin  I\left( {i,j} \right) \cap I\left( {i,k} \right) $. Clearly we have 
$c_p^i \neq c_p^j=c_p^k $. Hence $p\notin I(i,j)\cup I(i,k)$ and $p\in\overline{I\left( {i,j} \right) \cup I\left( {i,k} \right)}$.  It follows that $$I\left( {j,k} \right) \subseteq  ({I\left( {i,j} \right) \cap I\left( {i,k} \right)} ) \cup \overline{I\left( {i,j} \right) \cup I\left( {i,k} \right)}.$$ This completes the proof.\qed

Let $\mathcal{F}$ be a family of subsets of an $n$-set $E$. Let \begin{eqnarray}
&&l = \min \{|F|:F \in \mathcal{F}\}, \notag \\
&&u = \max \{|F|:F \in \mathcal{F}\}\notag
\end{eqnarray} be the minimum size and maximum size of subsets of ${\cal F}$.  For $r\ge u$ and $s\le l$, the families
\begin{eqnarray}
\nabla_r (\mathcal{F}) = \{B \subseteq E:|B|=r,\exists F \in \mathcal{F},F\subseteq B\},  \notag\\
\Delta_s(\mathcal{F}) = \{B \subseteq E:|B|=s,\exists F \in \mathcal{F},B\subseteq F\}  \notag
\end{eqnarray}
are called the {\it $r$-shade} and {\it $s$-shadow} of $\mathcal{F}$, respectively. When  $\mathcal{F}$ is a family of $k$-subsets, the ($k+1)$-shade and the ($k-1$)-shadow are simply written as $\nabla (\mathcal{F})$ or $\Delta(\mathcal{F})$.

{\lemma\label{interDelta} \cite{Katona64} If ${\cal A}$ is an intersecting family of $k$-subsets of an $n$-set, then $|\Delta {\cal A}|\ge|{\cal A}|$.}

\begin{lemma}\label{shadesha} \cite[Corollary 2.3.2]{EngelS}
Let $\mathcal{F}$ be a family of $k$-subsets of an $n$-set where $k<n$ and $n\geq 3$.

(1) If $k\geq \lceil\frac{n}{2}+1\rceil$, then $|\Delta \mathcal{F}|-|\mathcal{F}|\geq k-1\geq  \lceil\frac{n}{2} \rceil$.

(2) If $k\leq \lfloor\frac{n}{2}-1\rfloor$, then $|\nabla \mathcal{F}|-|\mathcal{F}|\geq n-k-1\geq  \lceil\frac{n}{2} \rceil$.
\end{lemma}

\begin{theorem}\label{u-l}
 Let $\mathcal{F}$ be a Sperner family over an $n$-set and let  $l\leq \frac{n}{2}\leq u$, where  $l$ and $u$ are the minimum size and the maximum size of subsets in ${\cal F}$, respectively.

(i) If $l=\lfloor {n\over 2}\rfloor$ and $u=\lceil{n\over 2}\rceil$, then $|{\cal F}|\le {n \choose \lfloor{\frac{n}{2}}\rfloor}.$

(ii) If $l<\lfloor {n\over 2}\rfloor$ and $u>\lceil{n\over 2}\rceil$, then \begin{eqnarray}
|\mathcal{F}| \leq
\begin{cases}
{n \choose {\frac{n}{2}}}-(u-l){n\over 2}-\lfloor{(u-l-1)^2\over 4}\rfloor,~~\mbox{if } n~\mbox{is even},\\
{n \choose {\frac{n-1}{2}}}-(u-l-1){n+1\over 2}-\lfloor{(u-l-2)^2\over 4}\rfloor,~~\mbox{if } n~\mbox{is odd}. \notag
\end{cases}
\end{eqnarray}

(iii) If $l<\lfloor {n\over 2}\rfloor$ and $u=\lceil{n\over 2}\rceil$, then $|{\cal F}|\le {n \choose \lfloor{\frac{n}{2}}\rfloor}-(\lfloor{n\over 2}\rfloor-l)\lceil{n\over 2}\rceil-{(\lfloor{n\over 2}\rfloor-l)(\lfloor{n\over 2}\rfloor-l-1)\over 2}.$

(iv) If $l=\lfloor {n\over 2}\rfloor$ and $u>\lceil{n\over 2}\rceil$, then $|{\cal F}|\le {n \choose \lfloor{\frac{n}{2}}\rfloor}-(u-\lceil{n\over 2}\rceil)\lceil{n\over 2}\rceil-{(u-\lceil{n\over 2}\rceil)(u-\lceil{n\over 2}\rceil-1)\over 2}.$
\end{theorem}

\proof  This is an adaption of {\rm \cite[Corollary 2.3.3]{EngelS}}.
The cases $n=1,2$ are trivial, hence we let $n\ge 3$. The statement (i) follows immediately from Sperner's Theorem.

In the case that $l<\lfloor {n\over 2}\rfloor$ and $u>\lceil{n\over 2}\rceil$, we proceed in two steps.

{\bf Step 1}: Replace ${\cal F}$ by ${\cal F}_1=({\cal F}\setminus{\cal G})\cup\nabla({\cal G})$ where ${\cal G}={\cal F}\cap{[n]\choose l}$. Because ${\cal F}$ is Sperner, we have that $({\cal F}\setminus{\cal G})\cap\nabla({\cal G})=\emptyset$ and that ${\cal F}_1$ is a Sperner family for which by Lemma \ref{shadesha} we have $$|{\cal F}_1|=|{\cal F}|-|{\cal G}|+|\nabla({\cal G})|\ge  |{\cal F}|+n-l-1.$$

If $l+1\le\lfloor{n\over 2}-1\rfloor$ then replace ${\cal F}_1$ by ${\cal F}_2=({\cal F}_1\setminus{\cal G}_1)\cup\nabla({\cal G}_1)$ where ${\cal G}_1={\cal F}_1\cap{[n]\choose l+1}$. After this we obtain a  Sperner family ${\cal F}_2$ for which  by Lemma \ref{shadesha} $$|{\cal F}_2|\ge |{\cal F}_1|+n-l-2\ge |{\cal F}|+(n-l-1)+(n-l-2).$$ Repeat doing like this until  we raise the minimum size of the subsets to $\lfloor{n\over 2}\rfloor$ and we obtain a  Sperner family ${\cal F}_{\lfloor{n\over 2}\rfloor-l}$ satisfying
\begin{align}\label{incr}|{\cal F}_{\lfloor{n\over 2}\rfloor-l}|\ge|{\cal F}|+(n-l-1)+(n-l-2)+\cdots+\lceil{n\over 2}\rceil.\end{align}

{\bf Step 2}:  We  begin to decrease the maximum size of subsets of ${\cal F}_{\lfloor{n\over 2}\rfloor-l}$.   Replace ${\cal F}_{\lfloor{n\over 2}\rfloor-l}$ by ${\cal F}_{\lfloor{n\over 2}\rfloor-l+1}=({\cal F}_{\lfloor{n\over 2}\rfloor-l}\setminus{\cal H})\cup\Delta({\cal H})$ where ${\cal H}={\cal F}_{\lfloor{n\over 2}\rfloor-l}\cap{[n]\choose u}$. Then we obtain a Sperner family ${\cal F}_{\lfloor{n\over 2}\rfloor-l+1}$ for which by Lemma \ref{shadesha} we have $$|{\cal F}_{\lfloor{n\over 2}\rfloor-l+1}|\ge |{\cal F}_{\lfloor{n\over 2}\rfloor-l}|+u-1.$$

If $u-1\ge\lceil{n\over 2}+1\rceil$ then replace ${\cal F}_{\lfloor{n\over 2}\rfloor-l+1}$ by ${\cal F}_{\lfloor{n\over 2}\rfloor-l+2}=({\cal F}_{\lfloor{n\over 2}\rfloor-l+1}\setminus{\cal H}_1)\cup\Delta({\cal H}_1)$ where ${\cal H}_1={\cal F}_{\lfloor{n\over 2}\rfloor-l+1}\cap{[n]\choose u-1}$. Similarly we have $$|{\cal F}_{\lfloor{n\over 2}\rfloor-l+2}|\ge |{\cal F}_{\lfloor{n\over 2}\rfloor-l+1}|+u-2\ge |{\cal F}_{\lfloor{n\over 2}\rfloor-l}|+(u-1)+(u-2).$$

Repeat this process until  we obtain a  Sperner family ${\cal F}_{\lfloor{n\over 2}\rfloor-l+u-\lceil{n\over 2}\rceil}$ with maximum size of the subsets being $\lceil{n\over 2}\rceil$ (and all sizes of the subsets being $\lfloor{n\over 2}\rfloor$ or $\lceil{n\over 2}\rceil$) and we have \begin{align}\label{decr1}{n \choose {\lfloor\frac{n}{2}}\rfloor}\ge|{\cal F}_{\lfloor{n\over 2}\rfloor-l+u-\lceil{n\over 2}\rceil}|\ge |{\cal F}_{\lfloor{n\over 2}\rfloor-l}|+(u-1)+(u-2)+\cdots+\lceil{n\over 2}\rceil.\end{align}
Combining inequality (\ref{decr1}) with (\ref{incr}) yields
 \begin{align}\label{decr}{n \choose {\lfloor\frac{n}{2}}\rfloor}\ge
|{\cal F}|+(n-l-1)+(n-l-2)+\cdots+\lceil{n\over 2}\rceil
 +(u-1)+(u-2)+\cdots+\lceil{n\over 2}\rceil.\end{align}

Whenever $n$ is even, 
 we further bound (\ref{decr}) by \begin{align*}{n \choose \frac{n}{2}}&\ge|{\cal F}|+(u-l){n\over 2}+(1+2+\cdots+(\lfloor{u-l\over 2}\rfloor-1))
+(1+2+\cdots+(\lceil{u-l\over 2}\rceil-1))
\\
&=|{\cal F}|+(u-l){n\over 2}+\lfloor{(u-l-1)^2\over 4}\rfloor.\end{align*}
%
Whenever $n$ is odd, similarly  we have \begin{align*}{n \choose {\frac{n-1}{2}}}&\ge|{\cal F}|+(u-l-1){n+1\over 2}+(1+2+\cdots+(\lfloor{u-l-1\over 2}\rfloor-1))\\
&\qquad\ \ +(1+2+\cdots+(\lceil{u-l-1\over 2}\rceil-1))\\
&=|{\cal F}|+(u-l-1){n+1\over 2}+\lfloor{(u-l-2)^2\over 4}\rfloor.\end{align*}
Then the statement (ii) follows immediately.

 Finally we consider (iii) and (iv). If $l<\lfloor {n\over 2}\rfloor$ and $u=\lceil{n\over 2}\rceil$,  then we only need to proceed Step 1 and the  Sperner family ${\cal F}_{\lfloor{n\over 2}\rfloor-l}$ consists of subsets of sizes $\lfloor{n\over 2}\rfloor$ and $\lceil{n\over 2}\rceil$. Thus from inequality (\ref{incr}) we have
\begin{align*} &{n\choose \lfloor{n\over 2}\rfloor}\ge|{\cal F}|+(n-l-1)+(n-l-2)+\cdots+\lceil{n\over 2}\rceil\\
&\quad\qquad=|{\cal F}|+(\lfloor{n\over 2}\rfloor-l)\lceil{n\over 2}\rceil+{(\lceil{n\over 2}\rceil-l)\rceil(\lceil{n\over 2}\rceil-l-1)\over 2},\end{align*}
proving (iii). Similarly we only proceed Step 2 and prove (iv) for the case that $l=\lfloor {n\over 2}\rfloor$ and $u>\lceil{n\over 2}\rceil$.
\qed


\begin{corollary}\label{i-shad}
 Let $\mathcal{F}$ be a Sperner family over an $n$-set. Let $l$ and $u$ be the minimum size and the maximum size of subsets in ${\cal F}$, respectively.
For $l\leq i\leq u$, define
\begin{eqnarray}
&&\mathcal{F}_{+}^i = \{B:~B \in \mathcal{F},|B|\geq i\}, \notag \\
&&\mathcal{F}_{-}^i = \{B:~B \in \mathcal{F},|B|\leq i\}.\notag
\end{eqnarray}

(i) If $l\leq \lfloor\frac{n}{2}\rfloor$, then $|\nabla_{i} (\mathcal{F}_{-}^i)|\geq |\mathcal{F}_{-}^i|+(i-l)(n-i)$ for any $l\leq i\leq \lfloor\frac{n}{2}\rfloor$.

(ii) If $u\geq \lceil\frac{n}{2}\rceil$, then $|\Delta_{i} (\mathcal{F}_{+}^i)|\geq|\mathcal{F}_{+}^i|+i(u-i)$ for any $\lceil\frac{n}{2}\rceil\leq i\leq u$.

\end{corollary}

\proof   Obviously the conclusion (i) holds  if $l= \lfloor\frac{n}{2}\rfloor$ and (ii) holds if $u= \lceil\frac{n}{2}\rceil$. So we let $l< \lfloor\frac{n}{2}\rfloor$ in (i) and let $u>\lceil\frac{n}{2}\rceil$ in (ii).

(i) Analogous to the  proof of inequality (\ref{incr}) of Theorem \ref{u-l}, we increase the minimum size of the subsets of ${\cal F}_{-}^i$ to $i$ step by step and then we have
\begin{align*}|\nabla_{i} (\mathcal{F}_{-}^i)|\ge|\mathcal{F}_{-}^i|+(n-l-1)+(n-l-2)+\cdots+(n-i)\geq |\mathcal{F}_{-}^i|+(i-l)(n-i).\end{align*}

(ii) Analogous to the  proof of inequality (\ref{decr1}) of Theorem \ref{u-l}, we decrease the maximum size of the subsets of ${\cal F}_{+}^i$ to $i$ step by step and then we have
$$|\Delta_{i} (\mathcal{F}_{+}^i)|\geq|\mathcal{F}_{+}^i|+(u-1)+(u-2)+\ldots+i\geq|\mathcal{F}_{+}^i|+i(u-i).$$ \qed

\section{Improved upper bounds}

Theorem \ref{FP&Sp} establishes the relationship between a 2-wFP code and non $2$-covering Sperner families generated by all codewords. Improved upper bounds on the size of a 2-wFP code will be developed in this section.

\subsection{Length even}

\begin{lemma}\label{even1}
 Let $n\ge 6$ be even and $\mathcal{F}$ a non $2$-covering Sperner family over an $n$-set.  Denote $l$ and $u$ to be the minimum size and the maximum size of subsets in ${\cal F}$, respectively.

(i) If $u\ge \frac{n}{2}+1$, then $|\mathcal{F}| \leq{n \choose \frac{n}{2}-1}-\frac{n}{2}. $

(ii) If $l\leq \frac{n}{2}-2$, then $|\mathcal{F}| \leq{n \choose \frac{n}{2}-1}-\frac{n}{2}-1.$

(iii) If $u=l=\frac{n}{2}$, then  $|\mathcal{F}| \leq\frac{1}{2} {n \choose {\frac{n}{2}}}.$

\end{lemma}

\proof Let $\mathcal{F} =\mathcal{A} \cup \mathcal{B}$ be a non $2$-covering Sperner family on $[n]$, where
\begin{eqnarray}
&\mathcal{A} = \{A:A \in \mathcal{F},|A|\geq \frac{n}{2}\},\notag\\
&\mathcal{B} = \{B:B \in \mathcal{F},|B|\leq \frac{n}{2}-1\}.\notag
\end{eqnarray}

(i) Let $u \geq \frac{n}{2}+1$. 
 By Corollary \ref{i-shad}, we have
  \begin{eqnarray}
|\Delta_{\frac{n}{2}}(\mathcal{A})| \geq  |\mathcal{A}|+(u-\frac{n}{2})\frac{n}{2} \geq  |\mathcal{A}|+ \frac{n}{2}. \notag
\end{eqnarray}
Denote ${\cal P}=\Delta_{\frac{n}{2}}(\mathcal{A})$. Because ${\cal F}$ is non 2-covering, we have that ${\cal P}$ is intersecting. Then  by Lemma \ref{interDelta} we have $$|\Delta(\mathcal{P})|\ge |{\cal P}|\ge |\mathcal{A}|+ \frac{n}{2}.$$

By Theorem \ref{symchains}, all  subsets of $[n]$ can be partitioned into ${n \choose {\frac{n}{2}}}$ mutually disjoint symmetric chains. If ${\cal B}\ne\emptyset$, then replace each $B\in{\cal B}$ with $B'$ in the same symmetric chain, where $B\subseteq B'$ and $|B'|= \frac{n}{2}-1$. Thus we produce from ${\cal B}$ a new Sperner family $\mathcal{B}{'}$ of $({n\over 2}-1)$-subsets. Note also that $\Delta(\mathcal{P})\cap{\cal B}'=\emptyset$  because ${\cal F}={\cal A}\cup{\cal B}$ is Sperner. (Let ${\cal B}'=\emptyset$ if ${\cal B}=\emptyset$.)
As a result,
\begin{eqnarray}
|\mathcal{B}|=|\mathcal{B}'|\leq {n \choose {\frac{n}{2}-1}}-|\Delta(\mathcal{P})|. \notag
\end{eqnarray}
It follows that
\begin{eqnarray}
|\mathcal{F}| = |\mathcal{A}| + |\mathcal{B}| \leq |\mathcal{A}|+ {n \choose {\frac{n}{2}-1}} -|\Delta(\mathcal{P})|\leq {n \choose {\frac{n}{2}-1}}- \frac{n}{2}.\notag
\end{eqnarray}

(ii) Let $l\leq \frac{n}{2}-2$. Apply Corollary \ref{i-shad} (i) to ${\cal B}$. It follows that
\begin{eqnarray}
|\nabla_{\frac{n}{2}-1}(\mathcal{B})| \geq |\mathcal{B}|+({n\over 2}-1-l)(n-({n\over 2}-1))\ge|\mathcal{B}|+\frac{n}{2}+1.\notag
\end{eqnarray}

In the decomposition of the power set of ${[n]}$ into symmetric chains, if ${\cal A}\ne\emptyset$, then replace each $A\in{\cal A}$ by $A'$ of the same symmetric chain where $|A'|={n\over 2}$ to obtain a new family ${\cal A}'$ of ${n\over 2}$-sets. Since ${\cal F}$ is non 2-covering, ${\cal A}'$ is intersecting and thus $|\Delta({\cal A}')|\ge|{\cal A}'|=|{\cal A}|$ by Lemma \ref{interDelta}. Furthermore, it is easy to see that $\Delta({\cal A}')$ and $\nabla_{\frac{n}{2}-1}(\mathcal{B})$ are  disjoint because ${\cal F}$ is Sperner.  (Let ${\cal A}'=\emptyset$ if ${\cal A}=\emptyset$.) It follows that
\begin{eqnarray}
|\mathcal{F}| =|\mathcal{A}|+|\mathcal{B}| \leq |\Delta(\mathcal{A}')|+|\nabla_{\frac{n}{2}-1}(\mathcal{B})|-\frac{n}{2}-1 \leq {n \choose {\frac{n}{2}-1}} -\frac{n}{2}-1.\notag
\end{eqnarray}

(iii) If $u=l=\frac{n}{2}$, then $|{\cal F}|\leq\frac{1}{2} {n \choose {\frac{n}{2}}}$ because ${\cal F}$ is  non 2-covering. \qed

 For a family ${\cal F}$ of subsets of $[n]$, we define its {\it complement} by $\overline{\cal F}=\{\overline{F}:F\in{\cal F}\}$.

{\theorem\label{even} Let $n$ be even and $n\geq 8$. Suppose that $\mathcal{C}$ is an $(n,m,q)$ $2$-wFP code. Then
$$m\leq{n \choose \frac{n}{2}-1}-\frac{n}{2}+1.$$}

\proof For $1\le i\le m$, let ${\cal X}_i$  be the  coincidence family generated by the codeword $c^i\in{\cal C}$. Then each ${\cal X}_i$ is a non $2$-covering Sperner family by Theorem \ref{FP&Sp}. Take a fixed $i\in[m]$ and let $l$ and $u$ be the minimum size and the maximum size of subsets in ${\cal X}_i$, respectively.

By Lemma \ref{even1}, if  $u\ge \frac{n}{2}+1$, then $|\mathcal{X}_i| \leq{n \choose \frac{n}{2}-1}-\frac{n}{2}$; if $l\leq \frac{n}{2}-2$, then $|\mathcal{X}_i| \leq{n \choose \frac{n}{2}-1}-\frac{n}{2}-1$; if $u=l=\frac{n}{2}$, then  $|\mathcal{X}_i| \leq\frac{1}{2} {n \choose {\frac{n}{2}}}= {n-1 \choose {\frac{n}{2}-1}}.$  It is easy to show that $m=|\mathcal{X}_i|+1\leq{n \choose \frac{n}{2}-1}-\frac{n}{2}+1$ for these cases. So, to prove the conclusion,  we only need to let $l\ge \frac{n}{2}-1$, $u\le \frac{n}{2}$, and $(u,l)\ne(\frac{n}{2},\frac{n}{2})$. Thus we  only need to consider two cases  $u=l= \frac{n}{2}-1$ and $(u,l)=(\frac{n}{2},\frac{n}{2}-1)$. Let \begin{eqnarray}
&\mathcal{A} = \{A:A \in \mathcal{X}_i,|A|= \frac{n}{2}\},\notag\\
&\mathcal{B} = \{A:A \in \mathcal{X}_i,|A|= \frac{n}{2}-1\}.\notag
\end{eqnarray}

\noindent\underline{Case 1}: Let $u=l= \frac{n}{2}-1$. Then ${\cal X}_i={\cal B}$. We evaluate the upper bound of $m$ by considering whether ${\cal B}$ is intersecting.

  If ${\cal B}$ is intersecting, then consider its complement $\overline{\mathcal{B}}=\{\overline{B}:B \in \mathcal{B}\}$.
For any $\overline{A},\overline{B} \in \overline{\mathcal{B}}$,
\begin{eqnarray}
|\overline{A}\cap \overline{B} | = n-|A\cup B |=n-(|A|+| B |-|A\cap B |)=n-(n-2-|A\cap B|)\geq 3. \notag
\end{eqnarray}
Consequently $\overline{\mathcal{B}}$ is a 3-intersecting Sperner family. By Theorem \ref{k-int}, we have
\begin{eqnarray}
m=|\mathcal{X}_{i}|+1 = |\overline{\mathcal{B}}|+1 \leq {n \choose {\frac{n}{2}-2}}+1. \notag
\end{eqnarray}

  If ${\cal B}$ is not intersecting, then there exist $B_{1},B_{2} \in \mathcal{B}$ such that $B_{1}\cap B_{2} =\emptyset$. Suppose that
$B_{1}=I(i,j)$ and
$B_{2}=I(i,k) $ where $i\ne j,k$ and $1\le j,k\le m$.
By Lemma {\ref{BH}}, $$I(j,k)\subseteq (B_1\cap B_2)\cup\overline{B_1\cup B_2}=\overline{B_1\cup B_2}$$ and hence $|I(j,k)|\leq 2$. If $|I(j,k)|=0$ then ${\cal X}_j$ is not Sperner, contradicting  Theorem \ref{FP&Sp}. So we have $|I(j,k)|=1,2$, meaning that ${\cal X}_j$ contains a set of size 1 or 2. Obviously if ${\cal X}_j$ contains a singleton, then by Proposition \ref{single} we have
$|\mathcal{X}_{j}| \leq {n-1 \choose {\frac{n}{2}-1}} +1 $ and hence $m=|\mathcal{X}_{j}| +1\le {n-1 \choose {\frac{n}{2}-1}} +2 $.
If $|I(j,k)|=2$, then the minimum size $l_j$ of elements of ${\cal X}_j$ satisfies  $l_j\le 2\le  \frac{n}{2}-2$ whenever $n\ge 8$. Hence $|\mathcal{X}_{j}| \leq{n \choose \frac{n}{2}-1}-\frac{n}{2}-1$ by Lemma \ref{even1}. Thus $m\le {n \choose \frac{n}{2}-1}-\frac{n}{2}$. Comparing the upper bounds of $m$ for $n\ge 8$ shows $m\le {n \choose \frac{n}{2}-1}-\frac{n}{2}$ in Case 1.

\noindent\underline{Case 2}: Let $u= \frac{n}{2}$ and $l=\frac{n}{2}-1$. Then ${\cal X}_i={\cal A}\cup{\cal B}$.

If ${\cal B}$ is not intersecting, then we have $m\le {n \choose \frac{n}{2}-1}-\frac{n}{2}$ by similar discussions in Case 1.

If ${\cal A}$ and ${\cal B}$ are not cross-intersecting, then there exist $A\in{\cal A}$ and $B \in \mathcal{B}$ such that $A\cap B =\emptyset$. Suppose that
$A=I(i,j)$ and
$B=I(i,k) $ where $i\ne j,k$ and $1\le j,k\le m$.
By Lemma {\ref{BH}}, $|I(j,k)|\leq|\overline{A\cup B}|= 1$. Since ${\cal X}_j$ is Sperner,  we have $|I(j,k)|=1$. Then by Proposition \ref{single} we have
$|\mathcal{X}_{j}| \leq {n-1 \choose {\frac{n}{2}-1}} +1 $ and hence $m\le {n-1 \choose {\frac{n}{2}-1}} +2 $.

What remains to bound $m$ is the subcase that  ${\cal B}$ is intersecting and that ${\cal A}$ and ${\cal B}$ are cross-intersecting.  By Lemma \ref{shadesha},
$|\nabla(\mathcal{B})| \geq |\mathcal{B}|+\frac{n}{2}.$ Note that ${\cal X}_i$ is Sperner and non 2-covering. As a result, we have that ${\cal A}\cap \nabla(\mathcal{B})=\emptyset$ and that ${\cal F}:={\cal A}\cup \nabla(\mathcal{B})$ is a non $2$-covering family of $\frac{n}{2}$-subsets. Hence \begin{eqnarray}
\frac{1}{2}{n \choose {\frac{n}{2}}}\ge|\mathcal{F}| = |\mathcal{A}\cup \nabla(\mathcal{B})| = |\mathcal{A}|+|\nabla(\mathcal{B})| \geq |\mathcal{A}|+|\mathcal{B}|+\frac{n}{2},\notag
\end{eqnarray}
yielding that \begin{eqnarray}
m=|\mathcal{X}_{i}|+1 = |\mathcal{A}|+|\mathcal{B}|+1 \leq \frac{1}{2}{n \choose {\frac{n}{2}}} - \frac{n}{2}+1={n-1 \choose {\frac{n}{2}}-1} - \frac{n}{2}+1.\notag
\end{eqnarray}
It follows that  $m\le {n-1 \choose {\frac{n}{2}-1}} +2 $ in Case 2.

To sum up it is immediate that for even $n\ge 8$ we have\begin{eqnarray}
m \leq {n \choose {\frac{n}{2}-1}}-\frac{n}{2}+1. \notag
\end{eqnarray} This completes the proof. \qed

\subsection{Length odd}

\begin{lemma}{\label{ood2}}
Let $n$ be odd and $n\geq 7$. Suppose that  $\mathcal{F}$ is a non $2$-covering Sperner family on $[n]$. If  $|F|\geq \frac{n+1}{2}$ for all $F\in \mathcal{F}$, then
\begin{eqnarray}
|\mathcal{F}|\le {n\choose {\frac{n+3}{2}}}. \notag
\end{eqnarray}
\end{lemma}

\proof Since $\mathcal{F}$ is a non $2$-covering Sperner family, for any $A,B\in{\cal F}$ we have
\begin{eqnarray}
n-1 \geq |A \cup B| = |A| +|B|-|A \cap B| \geq n+1- |A \cap B|,\notag
\end{eqnarray} implying $|A \cap B|\ge 2$. As a result
$\mathcal{F}$ is a 2-intersecting Sperner family. Then applying Theorem \ref{k-int} yields the conclusion. \qed

\begin{lemma}{\label{ood1}}
Let $n$ be odd and $n\geq 7$. Suppose that  $\mathcal{C}$ is  an  $(n,m,q)$ $2$-wFP code and ${\cal X}_i$ is the conincidence family generated by the codeword $c^i\in{\cal C}$. If  there  is $ i\in[ m]$ such that  $|F|\leq \frac{n-1}{2}$ for all $F\in \mathcal{X}_i$, then
\begin{eqnarray}
m\le {n\choose {\frac{n+3}{2}}}+1. \notag
\end{eqnarray}
\end{lemma}

\proof By Theorem \ref{FP&Sp}, ${\cal X}_i$ is a non $2$-covering Sperner family. Let
$\mathcal{X}_i={\mathcal{A}_0}\cup {\mathcal{A}_1}$, where
\begin{eqnarray}
&& {\mathcal{A}_0}=\{ A:A \in \mathcal{X}_i,|A| \leq \frac{{n - 3}}{2}\}, \notag\\
&&{\mathcal{A}_1}= \{ A:A \in \mathcal{X}_i,|A| = \frac{{n - 1}}{2}\}. \notag
\end{eqnarray}

If ${\mathcal{A}_1}$ is not intersecting, then there exist $I(i,j),I(i,k)\in{\cal A}_1$ ($1\le j,k\le m,i\ne j, k$) such that $I(i,j)\cap I(i,k)=\emptyset$. Hence by Lemma \ref{BH}  we have
$I(j,k)\subseteq\overline{I(i,j)\cup I(i,k)}$. Because ${\cal X}_i$ is Sperner, we know that $|I(j,k)|=1$. Apply Proposition \ref{single} to have
\begin{eqnarray}
m= |\mathcal{X}_j|+1\le {{n-1}\choose {\frac{n-1}{2}}}+2.\notag
\end{eqnarray}

If ${\mathcal{A}_1}$ is intersecting, then  let ${\mathcal F}=\overline{{\mathcal{A}_0}}\cup \overline{{\mathcal{A}_1}}$.
Similarly to the proof of Lemma \ref{ood2}, we can check that  ${\mathcal F}$ forms a 2-intersecting Sperner family. Hence by  Theorem \ref{k-int} we have \begin{eqnarray}
m=|\mathcal{X}_i|+1=|{\mathcal F}|+1\le {n\choose {\frac{n+3}{2}}}+1. \notag
\end{eqnarray}

Noting that $ {{n-1}\choose \frac{n-1}{2}}\le {n\choose \frac{n+3}{2}}$ if $n\ge 7$ yields the conclusion. \qed

\begin{theorem}{\label{odd}}
Let  $\mathcal{C}$ be  an  $(n,m,q)$ $2$-wFP code with $n$ odd and $n \geq 7$. Then
\begin{eqnarray*}
m \le \begin{cases}
{n \choose {\frac{n-1}{2}}}-{n^2-9\over 8}-\lfloor{(n-5)^2\over 64}\rfloor,~~\mbox{if } n\equiv 1\pmod 4,\\
{n \choose {\frac{n-1}{2}}}-{(n+1)^2-8\over 8}-\lfloor{(n-3)^2\over 64}\rfloor,~~\mbox{if } n\equiv 3\pmod 4.
\end{cases}
\end{eqnarray*}
\end{theorem}

\proof Let $\mathcal{C}$ be an $(n,m,q)$ $2$-wFP code and $\mathcal{X}_{i}$ be the non $2$-covering Sperner family generated by the codeword $c^{i}$. Denote
\begin{eqnarray}
&l_i = \min \{|A|:A \in \mathcal{X}_{i}\}, \notag \\
 &u_i= \max \{|A|:A \in \mathcal{X}_{i}\},\notag\\
&d_i=u_i-l_i.\notag
\end{eqnarray}

If there is $1\le i\le m$ such that $u_i\le{n-1\over 2}$ or $l_i\ge {n+1\over 2}$, then by Lemmas \ref{ood2} and \ref{ood1}
\begin{align}\label{two-extr}
m= |\mathcal{X}_i|+1\le {n\choose {\frac{n+3}{2}}}+1.
\end{align}

Next we let $l_i\le {n-1\over 2}$ and $u_i\ge{n+1\over 2}$ for all $1\le i\le m$.
Denote $d=\min\{d_i:1\le i\le m\}$ and assume w.l.o.g. $d_1=d$. 
In the case that $d\ge {n+1\over 2}$, by Theorem \ref{u-l}, it is easy to see that 
 \begin{align}\label{d>}m=|{\cal X}_1|+1\le {n\choose {\frac{n-1}{2}}}-{n^2-1\over 4}-\lfloor{(n-3)^2\over 16}\rfloor+1.\end{align}
In the following proof we let $d\le {n-1\over 2}$ and then bound $m$. Let ${{\cal X}_1} = {{\cal A}_1} \cup {{\cal A}_2} \cup {\cal B}$, where
\begin{eqnarray}
&{{\cal A}_1} = \{ A:A \in {{\cal X}_1},|A|\leq \frac{{n - 3}}{2}\},\notag\\
&{{\cal A}_2} = \{ A:A \in {{\cal X}_1},|A| = \frac{{n - 1}}{2}\}, \notag\\
&{\cal B} = \{ B:B \in {{\cal X}_1},|B| \ge \frac{{n + 1}}{2}\}. \notag
\end{eqnarray}
We consider three cases as follows.

\noindent\underline{Case 1}: Let ${\mathcal{A}_2}\ne \emptyset$ be not intersecting.  Then there exist $A,B\in{\cal A}_2$  such that $A\cap B=\emptyset$ and $|A\cup B|=n-1$. Suppose that $A=I(1,j)$ and $B=I(1,k)$. Then by Lemma \ref{BH} we have $I(j,k)\subseteq (A\cap B)\cup\overline{A\cup B}=\overline{A\cup B}$, meaning that ${\cal X}_j$ contains a singleton or $\emptyset$. Since ${\cal X}_j$ is Sperner, it contains a singleton and hence by Proposition \ref{single} we have \begin{align} \label{a2notint}
m=|{\cal X}_1|+1\le  {{n-1}\choose {\frac{n-1}{2}}}+2.
\end{align}

\noindent\underline{Case 2}: Let ${\cal A}_1\ne \emptyset$ and let ${\cal A}_1$ and ${\cal B}$ be not cross-intersecting. Then there exist  $A\in{\cal A}_1,B\in{\cal B}$ such that $A\cap B=\emptyset$.  Clearly
we have $|A|\ge{n+1\over 2}-d$ ($2\le d\le {n-1\over 2}$) and $|B|\ge{n+1\over 2}$. Suppose that $A=I(1,j)$ and $B=I(1,k)$. Then by Lemma \ref{BH} we have $I(j,k)\subseteq
\overline{A\cup B}$. 
It follows that $|I(j,k)|\le n-({n+1\over 2}-d+{n+1\over 2})=d-1$, meaning that ${\cal X}_j$ contains an $r$-subset where $1\le r\le d-1$. Hence we have $d_j\ge{n+1\over 2}-d+1$. Now we apply Theorem \ref{u-l} to bound  $m$. Let $d_0$ be an integer with $1\le d_0\le {n-3\over 2}$.

Whenever $d\ge d_0+1$,  by Theorem \ref{u-l}, noting ${(d_0-1)(d_0-2)\over 2}\ge\lfloor{(d_0-1)^2\over 4}\rfloor$, we have $$|{\cal X}_1|\le {n\choose {\frac{n-1}{2}}}-d_0{n+1\over 2}-\lfloor{(d_0-1)^2\over 4}\rfloor.$$

Whenever $d\le d_0$, we have $d_j\ge{n+1\over 2}-d_0+1$ and similarly we have
$$|{\cal X}_j|
\le {n\choose {\frac{n-1}{2}}}-({n+1\over 2}-d_0){n+1\over 2}-\lfloor{({n-1\over 2}-d_0)^2\over 4}\rfloor.$$

In order to get a better upper bound of $m=|{\cal X}_1|+1=|{\cal X}_j|+1$, we take $d_0={n-1\over 4}$ if $n\equiv 1$ (mod 4) and  $d_0={n+1\over 4}$ if $n\equiv 3$ (mod 4).
By simple reduction we have \begin{eqnarray}
m \leq
\begin{cases}
{n \choose {\frac{n-1}{2}}}-{n^2-1\over 8}-\lfloor{(n-5)^2\over 64}\rfloor+1,~~\mbox{if } n\equiv 1\pmod 4,\\
{n \choose {\frac{n-1}{2}}}-{(n+1)^2\over 8}-\lfloor{(n-3)^2\over 64}\rfloor+1,~~\mbox{if } n\equiv 3\pmod 4.
\end{cases}
\end{eqnarray}

\noindent\underline{Case 3}: Let ${\cal A}_2$ be intersecting if ${\cal A}_2\neq\emptyset$  and let ${\cal A}_1$  and ${\cal B}$ be cross-intersecting if ${\cal A}_1\neq\emptyset$. 
Define $${\cal F}={\overline{{\cal A}_1}}\cup{\overline{{\cal A}_2}}\cup{\cal B},$$ where we let $\overline{\emptyset}=\emptyset$ if necessary.  Then we claim that ${{\cal F}}$ is a 2-intersecting Sperner family.

Obviously each family of ${\overline{{\cal A}_1}},{\overline{{\cal A}_2}}$ and ${\cal B}$ is Sperner. If ${\cal F}$ is not Sperner, then one of the following three possibilities would happen: (a) there is $A_1 \in {{\cal A}_1}$ and $A_2 \in {{\cal A}_2}$ such that $\overline{A_2} \subseteq  \overline{A}_1$, yielding $A_1 \subseteq A_2$ and contradicting the fact that $\mathcal{X}_{1}$ is Sperner; (b) there is $A \in {{\cal A}_1}$ and $B \in \mathcal{B}$ such that $B \subseteq  \overline{A}$ or $ \overline{A} \subseteq B $, yielding $A\cap B=\emptyset$ or $A\cup B\supseteq A \cup\overline{A} =[n]$ and contradicting the fact that
${{\cal A}_1}$ and ${\cal B}$ are cross-intersecting or that $\mathcal{X}_{1}$ is non 2-covering;   (c) there is $A \in {{\cal A}_2}$ and $B \in {{\cal B}}$ such that $\overline{A} \subseteq  B$, yielding $A \cup B\supseteq  A\cup \overline{A}=[n]$ and contradicting the fact that $\mathcal{X}_{1}$ is non 2-covering. It follows that ${\cal F}$ is  Sperner.

Next we show that ${\cal F}$ is 2-intersecting. Noting the size of subsets in each family ${\overline{{\cal A}_1}},{\overline{{\cal A}_2}}$ and ${\cal B}$, we readily check that (a) $\overline {{\cal A}_1}$ is 3-intersecting, (b)
 $\overline {{\cal A}_2}$ is 2-intersecting because  ${\cal A}_2$ is intersecting, (c)  ${\cal B}$ is 2-intersecting as ${\cal B}$ is non 2-covering, (d)   ${\overline{{\cal A}_1}}$ and
${\overline{{\cal A}_2}}$  are cross-2-intersecting because $|A| + |B| \le n-2$ for  $A\in{\cal A}_1,B\in{\cal A}_2$, and (e)  ${\overline{{\cal A}_1}}\cup{\overline{{\cal A}_2}}$ and ${\cal B}$ are cross-2-intersecting  because for all $A\in{\cal A}_1\cup{\cal A}_2,B\in{\cal B}$ we have $A\cap B\ne A$ (${\cal X}_i$ Sperner) and
\begin{eqnarray}
|\overline {{A}}  \cap {{B}} | = |B|- |{A} \cap {B}|\ge|B|- (|{A}|-1)\ge {n+1\over 2}-{n-3\over 2} \ge 2. \notag
\end{eqnarray}

 Now that ${\cal F}$ is  a 2-intersecting Sperner family. So by Theorem \ref{k-int} we have \begin{align} \label{a2int}m=|{\cal X}_1|+1|=|{\cal F}|+1\le{n\choose {n+3\over 2}}+1.\end{align}

Comparing the upper bounds in (\ref{two-extr})-(\ref{a2int}) when $n\ge 7$ yields the conclusion. \qed

\section{Concluding remarks}

In the narrow-sense model, various types of fingerprinting codes were extensively studied by many researchers. However, as far as we know,  frameproof codes is the only type of fingerprinting codes that was ever studied  in the wide-sense model.
We tried the standard probabilistic method (as in \cite{Shangg1,StinsonWeiChen}) and approaches of  hypergraphs (as in \cite{Yangyt}) to achieve
lower bounds for wide-sense 2-frameproof codes. Unfortunately, we do not achieve a better lower bound than the bound of binary 2-FP codes in \cite [Theorem 4.1]{StinsonWeiChen}. Evaluating a good lower bound for 2-wFP code with general alphabet size  is a focus of future work.

The main result of the paper is an improvement on the known upper bounds for 2-wFP codes by applying techniques  on  non $2$-covering Sperner  families and intersecting families  in extremal set theory. The new bounds (as well as the previous ones by Panoui) do not relate with the alphabet size $q$. Not surprisingly, we learn from some small examples that $q$ usually affects the size of a $t$-wFP code. It is worthwhile to  examine how to improve the upper bounds by taking the alphabet size into consideration.


\medskip
\noindent{\bf Acknowledgements:} The authors would like to thank the anonymous referees
for their helpful comments and valuable suggestions, which have greatly improved the
presentation and quality of this paper.

\end{document}